\documentclass[12pt]{article}
\usepackage{amsfonts,amsmath,amssymb,epsfig,latexsym,rotating,mathrsfs}

\newtheorem{prop}{Proposition}[section]
\newtheorem{algor}[prop]{Algorithm} 
\newtheorem{thm}[prop]{Theorem}
\newtheorem{rem}[prop]{Remark}
\newtheorem{lemma}[prop]{Lemma}
\newtheorem{dfn}[prop]{Definition}

\newtheorem{baker}{Baker's Theorem (over $\Q$)}

\newtheorem{cor}[prop]{Corollary}

\newtheorem{conj}[prop]{Conjecture}
\newtheorem{not*}[prop]{Notation}

\newcommand{\N}{\mathbb{N}}
\newcommand{\Q}{\mathbb{Q}}
\newcommand{\R}{\mathbb{R}}
\newcommand{\C}{\mathbb{C}}

\newcommand{\Z}{\mathbb{Z}}
\newcommand{\Zn}{\Z^n}

\newcommand{\cA}{\mathcal{A}}
\newcommand{\sA}{\mathscr{A}}

\newcommand{\thth}{^{\text{\underline{th}}}}

\newcommand{\np}{\mathbf{NP}}

\newcommand{\sign}{\operatorname{sign}}

\newcommand{\qed}{$\blacksquare$}
\newcommand{\dia}{$\diamond$}
\newcommand{\eps}{\varepsilon}

\newcommand{\floor}[1]{\left\lfloor #1 \right\rfloor}

\begin{document}

\title{Trinomials and Deterministic 
Complexity Limits for Real Solving}

\maketitle

\author{Emma Boniface}\footnote{{\tt eboniface@berkeley.edu} . Partially supported by NSF REU grant 
DMS-1757872 and the Texas A\&{}M Mathematics Department.}
\author{Weixun Deng}\footnote{{\tt deng15521037237@tamu.edu} .  
Partially supported by NSF grant CCF-1900881.}
\author{J.\ Maurice Rojas}\footnote{{\tt rojas@tamu.edu} .  
Partially supported by NSF grant CCF-1900881.}

\begin{abstract}
We detail an algorithm that --- for all but a $\frac{1}{\Omega(\log(dH))}$ 
fraction of $f\!\in\!\Z[x]$ with exactly $3$ monomial terms, degree $d$, and 
all coefficients in $\{-H,\ldots, H\}$ --- produces an approximate root (in the 
sense of Smale) for each real root of $f$ in deterministic time $\log^{4+o(1)}(dH)$ 
in the classical Turing model. (Each approximate root is a rational with logarithmic height $O(\log(dH))$.) The best previous deterministic bit complexity 
bounds were exponential in $\log d$. We then relate this to Koiran's 
Trinomial Sign Problem (2017): Decide the sign of a degree $d$ trinomial $f\!\in\!\Z[x]$ with coefficients in $\{-H,\ldots,H\}$, at a point $r\!\in\!\Q$ of 
logarithmic height $\log H$, in (deterministic) time $\log^{O(1)}(dH)$. We show that 
Koiran's Trinomial Sign Problem admits a positive solution, at least for a 
fraction $1-\frac{1}{\Omega(\log(dH))}$ of the inputs $(f,r)$.  
\end{abstract}

\section{Introduction} 
The applications of solving systems of real polynomial equations 
permeate all of non-linear optimization, as well as numerous problems in 
engineering. As such, it is important to find the best possible speed-ups for real-solving. Furthermore, {\em structured} systems --- such as those with a fixed 
number of monomial terms or invariance with respect to a group action --- 
arise naturally in many computational geometric applications, and 
their computational complexity is closely related to a deeper understanding of circuit complexity 
(see, e.g., \cite{koiranrealtau}). 
So if we are to fully understand the complexity of solving sparse 
polynomial systems over the real numbers, then we should 
at least be able to settle the univariate case, e.g., classify when it is 
possible to separate and approximate roots in deterministic time polynomial 
in the input size. Independent of a complete classification, the underlying 
analytic estimates should give us a more fine-grained understanding of how 
randomization helps speed up real-solving for more general sparse polynomial 
equations. 

Recall that for any function $g$ analytic on $\R$, the 
corresponding {\em Newton endomorphism} is $N_g(z):=z-\frac{g(z)}{g'(z)}$, and the 
corresponding sequence of {\em Newton iterates} of a $z_0\!\in\!\R$
is the sequence $(z_i)^\infty_{i=0}$ where $z_{i+1}\!:=\!N_g(z_i)$ 
for all $i\!\geq\!0$. Given a trinomial 
$f(x)\!:=\!c_1+c_2x^{a_2}+c_3x^{a_3}\!\in\!
\Z[x]$ with $a_2\!<\!a_3\!=:\!d$ and all $c_i\!\in\!\{-H,\ldots,-1,1,\ldots,H\}$, 
we call $f$ {\em ill-conditioned} if and only if 
\begin{eqnarray}
\label{ineq:disc} 
\left|\left|\frac{c_2}{a_3}\right|\left|\frac{a_3-a_2}{c_1}\right|^{(a_3-a_2)/a_3}
\left|\frac{a_2}{c_3}\right|^{a_2/a_3}-1\right|<\frac{1}{\log(dH)}\\ 
\text{ and } f \text{ has no degenerate real roots.} \nonumber 
\end{eqnarray} 
We will see soon that Inequality (\ref{ineq:disc}) is the same as forcing the discriminant of $f$ to be near $0$ in an explicit way. Also, we'll see how we can check in time $\log^{2+o(1)}(dH)$ whether $f$ is ill-conditioned in the sense above. A peculiarity to observe that is that approximating real degenerate roots {\em is} also doable efficiently in our framework: Being ``near'' degeneracy --- not degeneracy itself --- is the remaining problem.  

We use $\#S$ for the cardinality of a set $S$. 
\begin{thm} {\em 
\label{thm:big} Following the notation above, assume $f$ is {\em not} ill-conditioned. Then we can find, in deterministic time $\log^{4+o(1)}(dH)$, a 
set $\left\{\frac{r_1}{s_1},\ldots,\frac{r_m}{s_m}\right
\}\!\subset\!\Q$ of cardinality $m\!=\!m(f)$ such that:\\   
\mbox{}\hspace{.5cm}1. For all $j$ we have $r_j\!\neq\!0 
\Longrightarrow  \log|r_j|,\log|s_j|=O(\log(dH))$.\\  
\mbox{}\hspace{.5cm}2. $z_0\!:=\!r_j/s_j\Longrightarrow f$ has a root $\zeta_j\!\in\!\R$ with sequence of Newton\\
\mbox{}\hspace{.8cm}iterates ($z_{i+1}\!:=\!N_{f'}(z_i)$ 
or $z_{i+1}\!:=\!N_f(z_i)$, according as $\zeta$ is\\
\mbox{}\hspace{.8cm}degenerate or not) satisfying $|z_i-\zeta_j|\!\leq\!(1/2)^{-2^{i-1}}|z_0-\zeta_j|$\\ 
\mbox{}\hspace{.8cm}for all $i\!\geq\!1$.\\  
\mbox{}\hspace{.5cm}3. $m\!=\!\#\{\zeta_1,\ldots,\zeta_m\}$ is the number of real roots of $f$.\\ 
In particular, if the exponents $a_2$ and $a_3$ are fixed (and distinct), then at most a 
fraction of $\frac{1}{\log(dH)}+\frac{1}{H}$ of the $(c_1,c_2,c_3)
\!\in\!\{-H,\ldots,H\}^3$ yield $f(x)\!=\!c_1+c_2x^{a_2}+c_3x^{a_3}$ that are  ill-conditioned. } 
\end{thm} 

\noindent 
We prove Theorem \ref{thm:big} in Section \ref{sec:trinosolr}, via 
Algorithm \ref{algor:trinosolr} there.  
We will call the convergence condition on $z_0$ above 
{\em being an approximate root (in the sense of Smale) 
with associated true root $\zeta_j$}. This type 
of convergence provides an efficient encoding of an approximation that 
can be quickly tuned to any desired accuracy. 
It is known (e.g., already for the special case of solving $x^2\!=\!c$) 
that one can not do much better, with respect to asymptotic arithmetic 
complexity, than Newton iteration \cite{bshouty}.  

Our complexity bound from Theorem \ref{thm:big} appears to be new, and complements earlier work on 
the arithmetic complexity of approximating \cite{rojasye,sag14} and counting 
\cite{brs,koiransep} real roots of trinomials. In particular, 
Theorem \ref{thm:big} nearly settles a question of Koiran from 
\cite{koiransep} on the bit complexity of solving trinomial equations over 
the reals. One should also observe that the best general bit complexity 
bounds for solving real univariate polynomials are super-linear in $d$ and 
work in terms of $\eps$-approximation, thus requiring an extra parameter 
depending on root separation (which is not known a priori): see, e.g., 
\cite{lickteig,pan}. 

\begin{rem} {\em Defining the {\em input size} of a univariate polynomial 
$f(x)\!:=\!\sum^t_{i=1} c_i x^{a_i}\!\in\!\Z[x]$ as 
$\sum^t_{i=1}\log((|c_i|+2)(|a_i|+2))$ we see that Theorem \ref{thm:big} 
implies that one can solve ``most'' real univariate\linebreak 
\scalebox{.95}[1]{trinomial equations 
in deterministic time polynomial in the input size. \dia}} 
\end{rem} 
\begin{rem} {\em 
Efficiently solving univariate $t$-nomial equations over $\R$ 
in the sense of Theorem \ref{thm:big} is easier for $t\!\leq\!2$: 
The case $t\!=\!1$ is clearly trivial (with $0$ the only possible root) 
while the case $t\!=\!2$ is implicit in work 
on computer arithmetic from the 1970s (see, e.g., \cite{borwein}). 
We review this case in Theorem \ref{thm:binor} of Section \ref{sub:bisep} below. \dia } 
\end{rem} 

Efficiently counting real roots for trinomials turns out to be equivalent 
to a special case of Baker's classic theorem on linear forms in logarithms 
\cite{baker, nesterenko}, and we review this equivalence in Lemma \ref{lemma:count} below. Our approach to {\em approximating} roots (in the sense of Smale) is to apply {\em $\cA$-hypergeometric 
functions} \cite{passare} (briefly 
reviewed in Section \ref{sub:hyper}) and a combination of earlier analytic estimates of Ye \cite{ye}, Rojas and Ye \cite{rojasye}, and Koiran \cite{koiransep} (see Sections \ref{sub:alpha} and  \ref{sub:tetrasep}). 

An important question Koiran posed near the end of his paper \cite{koiransep} is whether one can determine the sign of a trinomial evaluated at a rational number in (deterministic) time polynomial in the input size. (Determining the sign of a $t$-nomial at an {\em integer} turns out to be doable in deterministic polynomial-time for all $t$ \cite{cks}.) We obtain a partial positive answer to Koiran's question, thanks to an equivalence between solving and sign determination that holds for trinomials:
\begin{cor}
\label{cor:koiran} {\em 
Following the notation above, suppose $u,v\!\in\!\Z$ with $|u|,|v|\!\leq\!H$, and $f$ is not ill-conditioned. Then 
we can determine the sign of $f(u/v)$ in time $O(\log^{O(1)}(dH))$.} 
\end{cor} 
\begin{lemma} {\em 
\label{lemma:equiv} 
Koiran's Sign Problem has a positive solution if and only if finding approximate roots (in the sense of Smale) for trinomials is doable in deterministic time $\log^{O(1)}(dH)$.  } 
\end{lemma} 

\noindent 
We prove Corollary \ref{cor:koiran} and Lemma \ref{lemma:equiv} in Section \ref{sec:equiv}. Our 
use of $\cA$-hypergeometric series thus provides an alternative to how bisection is used to start higher-order numerical methods. In particular, our approach complements another approach to computing signs of trinomials at rational points of ``small'' height due to Gorav Jindal (write-up available at his blog). 

\subsection{The Root Separation Chasm at Four Terms} 
\label{sub:tetrasep} 
Unfortunately, there are obstructions to solving univariate polynomial equations over $\R$ in polynomial-time when there are too many monomial terms. Indeed, the underlying root spacing changes 
dramatically already at $4$ terms. 
\begin{thm} \cite{mig95,sag14,yuyupadic} 
{\em \label{thm:tetra}
Consider the family of tetranomials 
\begin{align*}
\displaystyle{f_d(x):=x^d - 4^h x^2 + 2^{h+2}x - 4 } 
\end{align*}
with $h\!\in\!\N$, $h\!\geq\!3$, and $d\!\in\!\left\{4,\ldots,\floor{e^h}
\right\}$ even. Let $H\!:=\!4^h$. Then $f_d$ 
has distinct roots $\zeta_1,\zeta_2\!\in\!\R$ with 
$\log|\zeta_1-\zeta_2|\!=\!-\Omega(d\log H)$. 
In particular, the coefficients of $f_d$ all lie 
in $\Z$ and have bit-length $O(\log H)$.} 
\end{thm}

\noindent
While this result goes back to work of Mignotte \cite{mig95}, we point 
out that in \cite{yuyupadic} a more general family of polynomials was 
derived, revealing that the same phenomenon of tightly-spaced roots for tetranomials 
occurs over {\em all} characteristic zero local fields, e.g., 
the roots of tetranomials in $\Q_p$ (for $p$ {\em any} prime) can be 
exponentially close as a function of the degree. One may conjecture that the basin of attraction, for Newton's Method applied to a real root of a tetranomial, can also be exponentially small, but so far only the analogous statement over $\Q_p$ is proved \cite[Rem.\ 4.1]{yuyupadic}.  

Tight spacing of real roots is thus partial evidence against being able to find 
approximate roots in the sense of Smale --- with ``small'' height, as in our main theorem for trinomials --- for tetranomials. Fortunately, for our setting, binomials and trinomials have well-spaced roots as a function of $d$ and $H$: \begin{thm} {\em (See \cite[Prop.\ 2.4]{yuyupadic} and \cite{koiransep}.) 
\label{thm:spacing} 
If $f\!\in\!\Z[x]$ is a degree $d$ univariate $t$-nomial, with 
coefficients in $\{-H,\ldots,H\}$, then any two distinct roots $\zeta_1,\zeta_2\!\in\!\C$ satisfy\\  
$\log|\zeta_1-\zeta_2|\!>\!-[\log(d)+\frac{1}{d}\log H]$ 
or $\log|\zeta_1-\zeta_2|\!=-O(\log^3(dH))$, according as 
$t$ is $2$ or $3$. \qed}  
\end{thm} 

One should also recall the following refined bound on the norms of nonzero roots of trinomials:  
\begin{lemma} 
\label{lemma:cauchy} 
Suppose $f(x)\!=\!c_1+c_2x^{a_2}+c_3x^{a_3}\!\in\!\C[x]\setminus\{0\}$ and $c_1c_2c_3\!\neq\!0$. Then any root $\zeta\!\in\!\C$ of $f$ must also satisfy\\  
$\frac{1}{2}\min\left\{\left|\frac{c_1}{c_2}\right|^{\frac{1}{a_2}},\left|\frac{c_1}{c_3}\right|^{\frac{1}
{a_3}}\right\}\!<\!|\zeta|\!<\!2
\max\left\{\left|\frac{c_2}{c_3}\right|^{\frac{1}{a_3-a_2}},\left|\frac{c_1}{c_3}\right|^{\frac{1}{a_3}}\right\}$. \qed 
\end{lemma} 

\noindent 
Such bounds had their genesis in work of Cauchy and Hadamard in the 19$\thth$ century, and have since been extended to several variables via tropical geometry: See, e.g., \cite{aknr,epr}. 

\subsection{Going Beyond Univariate Trinomials}
It is curious that the sign of an arbitrary $t$-nomial $f\!\in\!\Z[x]$ with degree $d$ and coefficients in $\{-H,\ldots,H\}$, at an {\em integer}  $r\!\in\!\{-H,\ldots,H\}$, can be computed in polynomial-time \cite{cks}, while the extension to rational $r$ is still an open question. It is conceivable that (but still unknown if) computing such sign evaluations at rational points can exhibit a leap in complexity for some family of tetranomials, akin to Theorem \ref{thm:tetra}.  

Let us call a polynomial in $\Z[x_1,\ldots,x_n]$ having exactly $t$ terms in its monomial term expansion an {\em $n$-variate $t$-nomial}. 
It is worth recalling that merely deciding the existence of roots 
over $\R$ for $n$-variate $(n+n^\eps)$-nomials (with $n\!\in\!\N$ and $\eps\!>\!0$ 
{\em arbitrary}) is $\np$-hard \cite{brs}. 

However, there is a different way to generalize univariate trinomial equations: They are the $n\!=\!1$ case of $n\times n$ {\em circuit systems}: Consider a system of equations $F\!:=\!(f_1,\ldots,f_n)\!\in\!\Z\!\left[x^{\pm 1}_1,\ldots,
x^{\pm 1}_n\right]$ where the exponent vectors of all the $f_i$ are contained in a set $A\!\subset\!\Zn$ of cardinality $n+2$, with 
$A$ {\em not} lying in any affine hyperplane. Such an $A$ is called a {\em circuit} (the terminology coming from combinatorics, instead of complexity theory), and such systems have been studied from the point of view of real solving and fewnomial theory since 2003 (if not earlier): See, e.g., \cite{lrw,bbs,brs,ckt}. In particular, it has been known at least since \cite{bbs} that solving such systems over $\R$ reduces mainly to finding the real roots of univariate rational functions of the form 
\begin{eqnarray} 
\label{eqn:galedual} 
g(u):=\prod^{n+1}_{i=1}\left(\gamma_{i,1}u+\gamma_{i,0}\right)^{b_i}-1  
\end{eqnarray} 
where $\gamma_{i,j}\!\in\!\Q$ and $b_i\!\in\!\Z$ for all $i,j$. Given any $u,v\!\in\!\Z$ with  
$\gcd(u,v)\!=\!1$, we define the {\em logarithmic height of 
$u/v$} to be $h(u/v)\!:=\!\max\{|u|,|v|\}$. (We also set $h(0)\!:=\!0$.) We pose the following conjecture: 
\begin{conj} 
\label{conj:solve} {\em 
Following the preceding notation, we can find approximate roots (in the sense of Smale) for all the real roots of (\ref{eqn:galedual}), in time 
polynomial in $\log^n(BH)$, where $B\!:=\!\max_i|b_i|$ and $\log H\!:=\!\max_{i,j} h(\gamma_{i,j})$. } 
\end{conj} 

\noindent 
Recently, it was shown that one can count the real roots of circuit systems in deterministic polynomial-time, for any fixed $n$ \cite{ckt}: The proof reduced to proving the simplification of Conjecture \ref{conj:solve} where one only asks for {\em the number} of real roots of $g$. 
This provides some slight evidence for Conjecture \ref{conj:solve}. More to the point, the framework from \cite{ckt} reveals that proving Conjecture \ref{conj:solve} would be the next step toward polynomial-time real-solving for circuit systems for $n\!>\!1$. Such speed-ups are currently known only for binomial systems so far \cite{ppr}, since real-solving for arbitrary $n\times n$ systems still has exponential-time worst-case complexity when $n$ is fixed (see, e.g., 
\cite{ls}).

\section{Background} \label{sec:back} 
\subsection{Approximating Logarithms and Roots of Binomials} 
\label{sub:bisep}
Counting real roots for the binomial $c_1+c_2x^d$ (with $c_1,c_2\!\in\!\Z$ and $d\!\in\!\N$) 
depends only on the signs of the $c_i$ and the parity of $d$: 
From the Intermediate Value Theorem, it easily follows that 
the preceding binomial has real roots if and only if 
[$c_1\!=\!0\!\neq\!c_2$, $c_1c_2\!<\!0$, or [$d$ is odd and $c_2\!\neq\!0$]]. Also, two nonzero real roots are possible if and only if [$c_1c_2\!<\!0$ and $d$ is even].   
So we now quickly review the bit complexity of finding a positive rational approximate 
root (in the sense of Smale) for $f(x)\!:=\!c_2x^d-c_1$, with $c_1,c_2,d\!\in\!\N$. (The case of negative roots obviously reduces to the case of positive roots by considering $f(-x)$.) 

\scalebox{.9}[1]{First note that $f$ must have a root in the open interval $\left(0,\max\left\{\frac{c_1}{c_2},1\right\}\right)$.}\linebreak 
So we can check 
the sign of $f$ at the midpoint of this interval and then reduce to either the 
left interval $\left(0,\frac{1}{2}\max\left\{\frac{c_1}{c_2},1\right\}\right)$, or the right interval  $\left(\frac{1}{2}\max\left\{\frac{c_1}{c_2},1\right\},\max\left\{\frac{c_1}{c_2},1\right\}\right)$, and proceed 
recursively, i.e., via the ancient technique of bisection. The signs can be computed efficiently by rapidly approximating 
$d\log x+\log(c_2/c_1)$, and other expressions of this form, to sufficiently many bits of accuracy. 

To see how to do this, we should first observe that logarithms of rational numbers 
can be approximated efficiently in the following sense: Recall that the
binary expansion of $\floor{2^{\ell-1-\floor{\log_2 x}}x}$ forms the
{\em $\ell$ most significant bits} of an $x\!\in\!\R_+$. (So knowing
the $\ell$ most significant bits of $x$ means that one knows
$x$ up to a multiple in the closed interval $\left[(1+2^{-\ell})^{-1},1+2^{-\ell})\right]$.)
\begin{thm}
\label{thm:log} {\em 
\cite[Sec.\ 5]{dan}
Given any positive $x\!\in\!\Q$ of logarithmic height $h$, and
$\ell\!\in\!\N$ with $\ell\!\geq\!h$, we can compute $\floor{\log_2 \max\{1,\log |x|\}}$, and the $\ell$ most
significant bits of $\log x$, in time $O(\ell\log^2\ell)$.} \qed
\end{thm}

\noindent 
The underlying technique ({\em AGM Iteration}) dates back to Gauss and was refined for
computer use in the 1970s by many researchers
(see, e.g., \cite{brent,salamin,borwein}).
We note that in the complexity bound above, we are applying
the recent $O(n\log n)$ algorithm of Harvey and van der Hoeven for multiplying
two $n$-bit integers \cite{harvey}. Should we use a more practical (but
asymptotically slower) integer multiplication algorithm then the time can
still be kept at $O\!\left(\ell^{1.585}\right)$ or lower.

The next fact we need is that only a moderate amount of accuracy is needed for 
Newton Iteration to converge quickly to a $d\thth$ root. 
\begin{lemma} {\em \cite{ye}  
\label{lemma:ye} 
Suppose $\zeta^d\!=\!c$ with $c\!\in\!\R_+$ and $d\!\in\!\N$. 
Then any $z\!\in\!\R_+$ satisfying $|z-\zeta|\!\leq\!\frac{2c^{1/d}}{d-1}$ is an 
approximate root of $x^d-c$ with associated true root $\zeta$. \qed } 
\end{lemma} 

The key to using fast logarithm computation to efficiently extract approximate $d\thth$ roots of rational numbers will then be knowing how roughly one can approximate the logarithms. An explicit estimate follows from a famous result of Baker, more recently refined by Matveev:  
\begin{baker} {\em (See \cite{baker} and \cite[Cor.\ 2.3]{matveev}.) 
Suppose $\alpha_i\!\in\!\Q\setminus\!\{0\}$ and $b_i\!\in\!\Z\setminus\!\{0\}$ for all 
$i\!\in\!\{1,\ldots,m\}$. Let  
$B\:=\!\max_i\{|b_1|,\ldots, |b_m |\}$,  
$\log\sA_i\!:=\!\max\{h(\alpha_i),|\log\alpha_i|, 0.16\}$, and $\Lambda\!:=\!\sum^m_{i=1}b_i\log \alpha_i$, where we fix any 
suitable branch of $\log$ a priori. Then 
$\Lambda\!\neq\!0\Longrightarrow\log|\Lambda|>
-1.4\cdot m^{4.5} 30^{m+3}(1+\log B)\prod^m_{i=1}
\log \sA_i$. \qed } 
\end{baker}

Combining Theorem \ref{thm:log}, Lemma \ref{lemma:ye}, and Baker's Theorem, 
we easily obtain the following result: 
\begin{thm} 
\label{thm:binor} {\em 
Suppose $f\!\in\!\Z[x]$ is a univariate binomial of degree $d$ with coefficients in $\{-H,\ldots,-1,1,
\ldots,H\}$. Then, in time $\log^{2+o(1)}(dH)$, 
we can count exactly how many real roots $f$ has and, for any nonzero real root 
$\zeta$ of $f$, find a $z_0\!\in\!\Q$, with $\zeta z_0\!>\!0$ and 
bit-length $O(\log(dH))$, that is an approximate root of $f$ in the sense of Smale. \qed}  
\end{thm} 

\noindent 
Theorem \ref{thm:binor} is most likely known to experts. In particular, an analogue for the arithmetic complexity of random binomial systems appears in \cite{ppr}.  

We now set the groundwork for extending the preceding theorem to the trinomial case. 

\subsection{Discriminants and $\alpha$-Theory for Trinomials} 
\label{sub:alpha}
There are three obstructions to extending the simple approach to binomials from last section to trinomials: (1) computing signs of trinomials at rational points is not known to be doable in polynomial-time, (2) counting roots requires the computation of the sign of a discriminant,
(3) we need explicit estimates on how close a rational $z$ must be to a real root $\zeta$ before $z$ can be used as an approximate root in the sense of Smale.

Circumventing Obstruction (1) is covered in the next section, so let us now review how to deal 
with Obstructions (2) and (3).

First recall the special case of the {\em $\cA$-discriminant} \cite{gkz94} for trinomials: 
\begin{dfn} 
{\em Given any $a_2,a_3\!\in\!\N$ with $\gcd(a_2,a_3)\!=\!1$ and $a_2\!<\!a_3$, we define\\ 
\mbox{}\hfill 
$\Delta_{\{0,a_2,a_3\}}(c_1,c_2,
c_3)\!:=\!a^{a_2}_2 
(a_3-a_2)^{a_3-a_2}(-c_2)^{a_3}
-a_3^{a_3}c^{a_3-a_2}_1c^{a_2}_3$,
\hfill\mbox{}\\ 
and abbreviate with 
$\Delta(f)\!:=\!\Delta_{(0,a_2,a_3)}(c_1,c_2,c_3)$ when\\  
$f(x)\!=\!c_1+c_2x^{a_2}+c_3x^{a_3}$. \dia } 
\end{dfn} 
\begin{rem} {\em 
By dividing out by a suitable monomial, the vanishing of $\Delta(f)$ is clearly equivalent to a monomial (with integer exponents) being $1$. Taking logarithms, we then see that we can decide the sign of any trinomial discriminant as above in time $\log^{2+o(1)}(dH)$, by combining Baker's Theorem with the fast logarithm approximation from Theorem \ref{thm:log}. Similarly, $f$ being ill-conditioned is equivalent to 
$\Delta(f)\!=\!O\!\left((1+\frac{1}{\log(dH)})^d\right)$, which can also be checked in time $\log^{2+o(1)}(dH)$ by approximating logarithms. 
\dia} 
\end{rem} 
\begin{lemma} \label{lemma:count} 
{\em Following the notation above, suppose\\ 
\mbox{}\hfill 
$f(x)\!=\!c_1+c_2x^{a_2}+c_3x^{a_3}\!\in\!\R[x]$ with 
$c_1c_2c_3\!\neq\!0$,\hfill\mbox{}\\  
and set $\sign(f)\!:=\!(\sign(c_1),\sign(c_2),\sign(c_3))\!\in\!\{\pm\}^3$. Then $f$ has...
\begin{enumerate}
\item{no positive roots if and only if   
[$\sign(f)\!\in\!\{(+,+,+),(-,-,-)\}$ or 
[$\sign(f)\!\in\!
\{(+,-,+),(-,+,-)\}$ and $\sign(\Delta(f))\!=\!\sign(c_2)$]].} 
\item{a unique positive root if and only if\\ \mbox{}  
[$\sign(f)\!\in\!\{(-,+,+),(-,-,+),(+,-,-),(+,+,-)$] or\\ \mbox{}  
[$\Delta(f)\!=\!0$ and 
$\sign(f)\!\in\!\{(+,-,+),(-,+,-)\}$]].}
\item{exactly two positive roots if and 
only if [$\sign(f)\!\in\!
\{(+,-,+)$,\linebreak 
$(-,+,-)\}$  
and $\sign(\Delta(f))\!=\!-\sign(c_2)$].} 
\item{$\left(-\frac{a_2c_2}{a_3c_3}\right)^{1/(a_3-a_2)}$ as a positive degenerate root (with no other positive root for $f$)  if and only if [$\Delta(f)\!=\!0$ and 
$\sign(f)\!\in\!\{(+,-,+),(-,+,-)\}$. \qed}  
\end{enumerate} }
\end{lemma} 

\noindent 
Lemma \ref{lemma:count} follows easily from Descartes' Rule of Signs (see, e.g., \cite{sl54}) and Assertion (4) (see, e.g., \cite{brs}). Since deciding the sign of $\Delta(f)$ is clearly 
reducible to deciding the sign of a linear combination of logarithms, Lemma \ref{lemma:count} combined with Baker's Theorem thus enables us to efficiently count the positive roots of trinomials (as already observed in \cite{brs}). 

So now we deal with the convergence of Newton's Method in 
the trinomial case. 
\begin{dfn} {\em 
For any analytic function $f : \R \longrightarrow \R$, let
$\gamma(f,x)\!:=\!\sup\limits_{k\geq 2} \left| \frac{f^{(k)}(x)}
{k!f'(x)}\right|^{\frac{1}{k-1}}$. \dia } 
\end{dfn}
\begin{rem} {\em 
\label{rem:gamma}
It is worth noting that $1/\gamma(f,x_0)$ is a lower bound for the
radius of convergence of the Taylor series of $f$ about $x_0$, so
$\gamma(f,x_0)$ is finite whenever $f'(x_0)\!
\neq\!0$ \cite[Prop.\ 6, Pg.\ 167]{bcss}. \dia } 
\end{rem}

\noindent 
A globalized variant, $\Gamma_f$, of $\gamma(f,x)$ will help us quantify how near $z\!\in\!\R_+$ must be to a positive root $\zeta$ of a trinomial for $z$ to be an approximate root in the sense of Smale with associated true root $\zeta$: 
\begin{dfn} \label{dfn:gamma}
{\em Consider $f(x)\!=\!c_1+c_2x^{a_2}+c_3x^{a_3}\!\in\!\R[x]$ with $0\!<\!a_2\!<\!a_3$, $c_3\!>\!0\!>c_2$, and $c_1\!\neq\!0$. 
Let $x_1$ be the unique positive root of the derivative $f'$.
\begin{enumerate} 
\item{If $f$ has two positive roots then let $x_2$ be the 
unique positive root of $f''$ (or set 
$x_2\!:=\!0$ should $f''$ not have a 
positive root). Then set 
$\Gamma_f\!:=\!\max\left\{\sup\limits_{x\in (0,x_2)} x\gamma(f,x),
\sup\limits_{x\in (x_2,\infty)}(x-x_1)\gamma(f,x)
\right\}$.} 
\item{If $c_1\!<\!0$ then we set 
$\Gamma_f\!:=\!\sup\limits_{x\in (0,x_2)} x\gamma(f,x)$. \dia} 
\end{enumerate}} 
\end{dfn} 

Lemma \ref{lemma:count} tells us that Cases (1) and (2) in our definition above are indeed disjoint, and Baker's Theorem (combined with Theorem \ref{thm:log}) tells us that we can efficiently distinguish Cases (1) and (2). Later, we will see some simple reductions implying that Cases (1) and (2) above are really the only cases we need to prove our main results. 
\begin{thm} {\em 
(See \cite[Thm.\ 2]{ye} and \cite[Thm.\ 5]{rojasye}.)   
\label{thm:alpha} 
Following the notation and assumptions of Definition \ref{dfn:gamma}, set $d\!:=\!a_3$ and 
suppose $z,\zeta\!\in\!\R_+$ with $f(\zeta)\!=\!0$ and $d\!\geq\!3$. Also let $x_2$ 
be the unique positive root of $f''$ (or set $x_2\!:=\!0$ should $f''$ not have a 
positive root). Then: 
\begin{enumerate} \setcounter{enumi}{-1}
\item{$-f$ is convex on $(0,x_2)$ and $f$ is convex on $(x_2,\infty)$.} 
\item{If $f$ is monotonically decreasing in $[z,\zeta]$, $x_2\not\in [z,\zeta]$, and $\zeta\!\in\!\left[z,
\left(1+\frac{1}{8\Gamma_f}\right)z\right]$,  
then $z$ is an approximate root of $f$. }
\item{If $f$ is monotonically increasing in $[\zeta,z]$, $x_2\not\in [x,\zeta]$, and $\zeta\!\in\!\left[
\left(1-\frac{1}{8\Gamma_f}\right)z,z
\right]$, then $z$ is an approximate root of $f$.}
\item{$\frac{d-1}{2}\!\leq\!\Gamma_f\!\leq\!\frac{(d-1)(d-2)}{2}$.}  
\end{enumerate} 
In particular, Assertions (1) and (2), combined with Assertion (3), imply $|z-\zeta|\!\leq\!\frac{\zeta}{4(d-1)(d-2)}$. \qed
} 
\end{thm}

\bigskip 
\subsection{$\cA$-Hypergeometric Functions} \label{sub:hyper} 
Let us now consider the positive roots of  $1-cx^m+x^n$ and\linebreak 
$-1-cx^m+x^n$ as a function 
of $c\!\in\!\R$, when $0\!<\!m\!<\!n$ and $\gcd(m,n)\!=\!1$, from the point of view 
of $\cA$-hypergeometric series \cite{gkz94}. These series
date back to 1757 work of Johann Lambert for the special case
$m\!=\!1$. Many authors have since extended these series
in various directions. Passare and Tsikh's paper \cite{passare} 
is the most relevant for our development here. 
The union of the domains of convergence of these series will turn out to be 
$\C\setminus\!\{r_{m,n}\}$ where 
$r_{m,n}:=\frac{n}{m^{\frac{m}{n}}(n-m)^\frac{n-m}{n}}$: This quantity, easily checked to be strictly greater than $1$, is closely related to the trinomial discriminant via 
$\Delta(\pm 1-cx^m+x^n)\!=\!0 \Longrightarrow |c|=r_{m,n}$. 

There will be two series for $|c|\!>\!r_{m,n}$, and one more series 
for $|c|\!<\!r_{m,n}$, that we will focus on. These two convergence domains are in fact related to triangulations of the point set $\{0,m,n\}$ via the {\em Archimedean Newton polygon} \cite{aknr,gkz94}, but we will not elaborate further here on this combinatorial aspect. The reader should be aware that these series in fact yield all complex roots, upon inserting suitable roots of unity in the series formulae. So our choices of sign are targeted toward numerical bounds for positive roots. 

\bigskip
\bigskip
\noindent
\mbox{}\hfill
{\large {\sc $|c|\!>\!r_{m,n}\Longrightarrow$ finest lower hull \
\epsfig{file=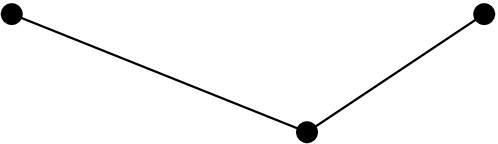,height=.15in}}}\hfill\mbox{}

\noindent
The first series is:\\  
\mbox{}\hfill
\fbox{$\displaystyle{x_\mathrm{low}(c)= \frac{1}{c^{1/m}}
\left[1+\sum^\infty_{k=1}\left(\frac{1}{km^k}\cdot \prod^{k-1}_{j=1}
\frac{1+kn-jm}{j}\right)\left(\frac{1}{c^{n/m}}
\right)^k\right]}$}
\hfill\mbox{}\\
(yielding a positive root of $1-cx^m+x^n$ with norm within a factor of $2$ of $c^{-1/m}$, thanks to Lemma \ref{lemma:cauchy}),
while the second is:\\
\mbox{}\hfill
\fbox{\scalebox{.8}[1]{$\displaystyle{x_\mathrm{hi}(c)=c^{\frac{1}{n-m}}\left[1- 
\sum^\infty_{k=1} \left(\frac{1}{k(n-m)^k}\cdot \prod^{k-1}_{j=1}
\frac{km+j(n-m)-1}{j}\right)\left(\frac{1}{c^{n/(n-m)}}
\right)^k \right]}$}}\hfill\mbox{}\\
(yielding a positive root of $1-cx^m+x^n$ with norm within a factor of $2$ of $c^{1/(n-m)}$, thanks to 
Lemma \ref{lemma:cauchy}, and distinct from the previous root). 

\bigskip
\bigskip 
\noindent
\mbox{}\hfill
{\large {\sc $|c|\!<\!r_{m,n}\Longrightarrow$ coarsest lower hull \
\raisebox{-.3cm}{\epsfig{file=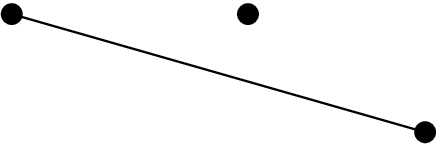,height=.15in,angle=16}}}}
\hfill\mbox{}

\vspace{-.3cm} 
\noindent 
Our final series yields a root of 
$-1-cx^m+x^n$ in the interval 
$\left(\frac{1}{2},2\right)$, via Descartes' Rule of Signs and Lemma \ref{lemma:cauchy}:\\
\mbox{}\hfill \fbox{$\displaystyle{x_\mathrm{mid}(c)=1+
\sum^\infty_{k=1}\left(\frac{1}{kn^k}\cdot \prod^{k-1}_{j=1} 
\frac{1+km-jn}{j} \right) c^k}$.}
\hfill\mbox{}

A key fact about these $\cA$-hypergeometric 
series is that their tails decay quickly enough 
for us to use their truncations (sometimes efficiently) as start points for Newton iteration.  

\vspace{-.1cm} 
\begin{lemma}
\label{lemma:tail}
{\em Suppose $f(x)\!=\!c_1+c_2x^{a_2}+c_3x^{a_3}\!\in\!\Z[x]$, with $\sigma(f)\!\in\!\{(-,-,+),(+,-,+)\}$ and all coefficients in $\{-H,\ldots,-1,1,\ldots,H\}$, is not ill-conditioned. Set $c\!:=\!\left|\frac{c_2}{
c^{(n-m)/n}_1 c^{m/n}_3}\right|$. Also 
let $x^{(\ell)}_\mathrm{low}$, 
$x^{(\ell)}_\mathrm{hi}$, and  
$x^{(\ell)}_\mathrm{mid}$ denote the 
truncation of the corresponding series to its 
$\ell\thth$ term. 
Then $|c_1/c_3|^{1/a_3}x^{(\ell)}_\mathrm{low}(c)$ (resp.\ $|c_1/c_3|^{1/a_3}
x^{(\ell)}_\mathrm{hi}(c)$,  $|c_1/c_3|^{1/a_3}x^{(\ell)}_\mathrm{mid}(c)$) is an approximate root of $f$ with associated true root $|c_1/c_3|^{1/a_3}x_\mathrm{low}(c)$ (resp.\ 
$|c_1/c_3|^{1/a_3}x_\mathrm{hi}(c)$, $|c_1/c_3|^{1/a_3}x^{(\ell)}_\mathrm{mid}(c)$) if $\ell=O(\log^2(dH))$. }
\end{lemma} 

\noindent 
{\bf Proof:} Let $m\!:=\!a_2$ and 
$n\!:=\!a_3$. Observe that if we set $\alpha\!:=\!|1/c_1|$,  $\beta\!:=\!|c_1/c_3|^{1/n}$, and 
$g(x)\!:=\!\alpha f(\beta x)$, then $g(x)\!=\!\pm 1 - cx^m+x^n$, with the sign being exactly $\sign(c_1)$. So we'll work with $g$ henceforth. Now note that\\
$\displaystyle{\prod_{j=1}^{k-1} \frac{1+kn-jm}{j}= \exp\left(\sum_{j=1}^{k-1} \log (1+kn-jm) - \log(j) \right)}$
\begin{align*}
&\leq \exp\left(\log(1+kn-1\cdot m) + \int_{1}^{k-1} \log (1+kn-jm) - \log(j) \ dj \right)\\
&= \left( \frac{1-m+kn}{1+m-km+kn}\right)^{\frac{1+kn}{m}} \left(\frac{1+m-km+kn}{k-1}\right)^{k-1}\\
&\leq k \left(m \ r_{m,n}^{\frac{n}{m}}\right)^k. 
\end{align*}
So $\displaystyle{\left|x_\mathrm{low}(c)-x^{(\ell)}_\mathrm{low}(c)\right|}$
\begin{align*}
&=\left|\frac{-1}{c^{1/m}}\sum_{k=\ell+1}^{\infty}\left(\frac{(-1)^{nk}}{km^k} \cdot \prod_{j=1}^{k-1} \frac{1+kn-jm}{j}\right)\left(\frac{1}{c^{n/m}}\right)^k\right| \\
&\leq \frac{1}{\left|c\right|^{1/m}}\sum_{k=\ell+1}^{\infty} \left(\frac{r_{m,n}^{n/m}}{|c|^{n/m}}\right)^k\\\
&= \frac{1}{\left|c\right|^{1/m}} \cdot \frac{1}{1-\left(r_{m,n}/\left|c\right|\right)^{\frac{n}{m}}}\cdot \left(\frac{r_{m,n}}{\left|c\right|}\right)^{\frac{n}{m}(\ell+1)}.\\
\end{align*}
Assume $c\!>\!r_{m,n}$ (and recall $r_{m,n}\!>\!1$). By Theorem \ref{thm:alpha}, when $n\!\geq\!3$, it suffices to find $\ell$ such that 
$\left|x_\mathrm{low}(c)-x^{(\ell)}_\mathrm{low}(c)\right|\leq \frac{x_\mathrm{low}(c)}{4(n-1)(n-2)}$. When $n\!=\!2$, we can in fact complete the square, reduce to the binomial case, and then Lemma \ref{lemma:ye} tells us that it suffices to find $\ell$ such that
$\left|x_\mathrm{low}(c)-x^{(\ell)}_\mathrm{low}(c)\right|\leq 2x_\mathrm{low}(c)$. Also, Lemma \ref{lemma:cauchy} tells us that\linebreak  $x_\mathrm{low}(c)\!\in\!\left(\frac{1}{2c^{1/m}},
\frac{2}{c^{1/m}}\right)$. 
In other words, it suffices to enforce\linebreak  
$\left|x_\mathrm{low}(c)-x^{(\ell)}_\mathrm{low}(c)\right|\leq \frac{1}{8(n-1)^2c^{1/m}}$, which (thanks to our last tail bound involving $r_{m,n}/|c|$) is  implied by:\\ 
\scalebox{1}[1]{$(\ell+1)\log\frac{\left|c\right|}{r_{m,n}}\geq\frac{m}{n}\log(8(n-1)^2)-\frac{m}{n}\log\left(1-\left(\frac{r_{m,n}}{\left|c\right|}\right)^{\frac{n}{m}}\right) 
$.}

\noindent 
Since $f$ is not ill-conditioned we have $\frac{\left|c\right|}{r_{m,n}}\geq 1+\frac{1}{\log(nH)}$ and thus 
\begin{align*}
-\frac{m}{n}\log\left(1-\left(\frac{r_{m,n}}{\left|c\right|}\right)^{\frac{n}{m}}\right)&\leq -\frac{m}{n}\log\left(1-\left(\frac{1}{1+\frac{1}{\log(nH)}}\right)^{\frac{n}{m}}\right)\\
&\leq \log\left(1+\frac{1}{\log(nH)}\right)
-\frac{m}{n}\log\frac{1}{\log(nH)}\\
&\leq \log(nH)
\end{align*}
Also,
$\log\frac{\left|c\right|}{r_{m,n}}\geq \log\left(1+\frac{1}{\log(nH)}\right)\geq \frac{\log2}{\log(nH)}$ 
by the inequality $\log(x)\geq (\log 2)(x-1)$ if $1\leq x\leq 2$. Therefore, if \[ \ell\geq\frac{\log(nH)[\log(8(n-1)^2)
+\log(nH)]}{\log 2},\] 
we then have $(\ell+1)\log\frac{\left|c\right|}{r_{m,n}}\!>\!\log(8(n-1)^2)+\log(nH)$
\begin{align*}
&\geq \frac{m}{n}\log(8(n-1)^2)-\frac{m}{n}\log\left(1-\left(\frac{r_{m,n}}{\left|c\right|}\right)^{\frac{n}{m}}\right).&
\end{align*}

The proof for $x^{(\ell)}_\mathrm{hi}(c)$ is similar, just using
\[
\prod^{k-1}_{j=1}\frac{1+kn-jm}{j}\leq k \left((n-m) \ r_{m,n}^{\frac{n}{n-m}}\right)^k
\]
and
$\left|x_\mathrm{hi}(c)-x^{(\ell)}_\mathrm{hi}(c)\right|\leq 
\frac{c^{\frac{1}{n-m}}}{1-\left(r_{m,n}/\left|c\right|\right)^{\frac{n-m}{n}}}\cdot \left(\frac{r_{m,n}}{\left|c\right|}\right)^{\left(\frac{n-m}{n}\right)(\ell+1)}$ instead.

The proof for $x^{(\ell)}_\mathrm{mid}(c)$ 
also follows similarly, assuming\linebreak 
$0\!<\!|c|\!<\!r_{m,n}$ instead, and using
\[\prod_{j=1}^{k-1} \frac{1+km-jn}{j} \leq k \left(\frac{n}{r_{m,n}}\right)^k
\] 
and $\left|x_{\mathrm{mid}}(c)-x_{\mathrm{mid}}\right| \leq \frac{1}{1-\left(c/r_{m,n}\right)} \left(\frac{c}{r_{m,n}}\right)^{\ell +1}$ instead.
\qed 

\section{Solving Trinomial Equations over $\R$} \label{sec:trinosolr} 
\scalebox{.9}[.9]{\fbox{\mbox{}\hspace{.3cm}\vbox{
\begin{algor} {\em
\label{algor:trinosolr}
{\bf (Solving Trinomial Equations Over $\pmb{\R_+}$)}
\mbox{}\\
{\bf Input.} $c_1,c_2,c_3,a_2,a_3\!\in\!\Z\setminus\{0\}$ with 
$|c_i|\!\leq\!H$ for all $i$ and\\  $1\!\leq\!a_2\!<\!a_3\!=:\!d$.\\
{\bf Output.} $z_1,\ldots,z_m\!\in\!\Q_+$ with logarithmic
height $O\!\left(\log(dH)\right)$ such that $m$ ($\leq\!2$) is the number of roots of $f(x):=c_1+c_2x^{a_2}+c_3x^{a_3}$ in $\R_+$, $z_j$ is an approximate root of
$f$ with associated true root $\zeta_j\!\in\!\R_+$ for all $j$, 
and the $\zeta_j$ are pair-wise distinct. \\
{\bf Description.}  
\begin{enumerate} 
\setcounter{enumi}{-1}
\item{Let {\tt xflip}$:=1$, $c\!:=\!\left|\frac{c_2}{c^{(a_3-a_2)/a_3}_1c^{a_2/a_3}_3}\right|$, 
$\beta\!:=\!|c_1/c_3|^{1/a_3}$,\linebreak and let 
$c'\!\in\!\Q_+$ (resp.\ $\beta'$) an approximation to $c$ (resp.\ $\beta$) within distance $\frac{c}{96(n-1)^2}$ computed via Theorem \ref{thm:binor}, and\linebreak 
$\ell\!:=\!\frac{\log(a_3H)
[\log(24(a_3-1)^2)+\log(a_3H)]}{\log 2}$.}
\item{If $\sigma(f)\!\in\!\{(+,+,+),(-,-,-)\}$ then {\tt output\\  ``Your $f$ has no positive roots.''}  
and {\tt STOP}.} 
\item{Replacing $f(x)$ by $\pm f(x)$ or  
$\pm x^{a_3}f(1/x)$ (and setting\linebreak 
{\tt xflip}$:=-1$) as necessary, 
reduce to the special case $c_3\!>\!0\!>\!c_2$.}  
\item{If $\Delta(f)\!=\!0$ then, using 
Theorem \ref{thm:binor}, let $z_1^{\mathtt{xflip}}$ be a rational approximation to   $\left(-\frac{a_2c_2}{a_3c_3}\right)^{1/(a_3-a_2)}$ of logarithmic height $O(\log(dH))$, {\tt output  
``$z_1$ is your only positive}\linebreak 
{\tt approximate root.''}  
and {\tt STOP}.}
\item{If $\Delta(f)\!<\!0$ then {\tt output}\\ 
{\tt ``Your trinomial has no positive roots.''} 
and {\tt STOP}.}
\item{If $c_1\!<\!0$ then let  $z^{\mathtt{xflip}}_1\!:=\!\beta' x^{(\ell)}_{\mathrm{mid}}(c')$ and {\tt output ``$z_1$ is your only positive approximate root.''} and 
{\tt STOP}.} 
\item{Let $z^{\mathtt{xflip}}_1\!:=\!\beta' x^{(\ell)}_{\mathrm{low}}(c')$, 
$z^{\mathtt{xflip}}_2\!:=\!\beta' x^{(\ell)}_{\mathrm{hi}}(c')$, and 
{\tt output ``$z_1$ and $z_2$ are your only positive approximate roots.''} and {\tt STOP}. } 
\end{enumerate} } 
\end{algor}}}} 

\medskip 
\noindent 
{\bf Proof of Theorem \ref{thm:big}:} 
We make one arithmetic reduction first: By computing $\delta\!:=\!\gcd(a_2,a_3)$ first, replacing $(a_2,a_3)$ with $(a_2,a_3)/\delta$, and solving the resulting trinomial $\bar{f}$, we can solve $f$ over $\R$ by taking the $\delta\thth$ root of all the real roots of $\bar{f}$ if $\delta$ is odd. (If $\delta$ is even, then we only take $\delta\thth$ roots of the positive roots of $\bar{f}$.) The underlying computation of $\delta\thth$ roots is done via rational approximations with precision $\frac{1}{96H(a_3-1)^2}$ via Theorem \ref{thm:binor}, possibly at the expense of a few extra Newton Iterations of neglible cost. The precision guarantees that the resulting approximations are indeed approximate roots in the sense of Smale for $f$. The computation of $\gcd(a_2,a_3)$ takes time $O(\log(d)(\log(\log d))^2)$ via the {\em Half-GCD Method} \cite{vzgbook}, so this reduction to the 
case $\gcd(a_2,a_3)\!=\!1$ has negligible 
complexity. 

Assuming Algorithm \ref{algor:trinosolr} is correct and runs within the stated time bound, our theorem then follows directly by applying Algorithm \ref{algor:trinosolr} to $f(x)$ and $f(-x)$. 
So it suffices to prove correctness, and analyze the complexity, of Algorithm \ref{algor:trinosolr}. 

\medskip 
\noindent 
{\bf Correctness:} This follows directly from Lemma \ref{lemma:count}, Theorem \ref{thm:alpha}, and Lemma \ref{lemma:tail}. In particular, the constants in our algorithm are chosen so that multiplicative error in the underlying radicals and evaluations of our series combine so that $\beta'x^{(\ell)}_\mathrm{mid}(c')$, 
$\beta'x^{(\ell)}_\mathrm{low}(c')$, and 
$\beta'x^{(\ell)}_\mathrm{hi}(c')$ are indeed approximate roots in the sense of Smale. 

\medskip 
\noindent 
{\bf Complexity Analysis:} The logarithmic height bounds on our approximate roots follow directly from construction, and our time bound follows easily upon observing that the truncated series we evaluate involve only $O(\log^2(dH))$ many  terms, and each term is an easily computable rational multiple of the previous one. In particular, to get $O(\log(dH))$ bits of accuracy it suffices to compute the leading $O(\log(dH))$ bits of each term (with a suitable increase of the second $O$-constant).  

Our final assertion on the fraction of trinomials that are ill-conditioned can be obtained as follows: We want to find the cardinality of the set:\\ 
\scalebox{.94}[1]{$\displaystyle{
\left\{\left(c_1,c_2,c_3\right)\in \left\{-H,\ldots,H\right\}^3\; \left|\; 
0\!<\!\frac{\left|c_2\right|}{\left|c_1\right|^{\frac{n-m}{n}}\left|c_3\right|^{\frac{m}{n}}r_{m,n}}-1\!<
\!\frac{1}{\log \left(dH\right)}\right.\right\}.}$} 

Fix $\left(c_1,c_3\right)\in \left\{-H,\ldots,H\right\}^2 $. Then $c_2$ satisfies
\[
r_{m,n}\left|c_1\right|^{\frac{n-m}{n}}\left|c_3\right|^{\frac{m}{n}}<\left|c_2\right|<\left(1+\frac{1}{\log \left(dH\right)}\right)r_{m,n}\left|c_1\right|^{\frac{n-m}{n}}\left|c_3\right|^{\frac{m}{n}}
\]
Since $\left|c_2\right|$ is an integer, the number of $c_2$ satisfying the last inequality is no more than
\begin{align*}
\ \ &\sum_{
\left(c_1,c_3\right)\in \left\{-H,\ldots,H\right\}^2
\atop
r_{m,n}\left|c_1\right|^{\frac{n-m}{n}}\left|c_3\right|^{\frac{m}{n}}\leq H}
2\left \lceil \frac{1}{\log \left(dH\right)}r_{m,n}\left|c_1\right|^{\frac{n-m}{n}}\left|c_3\right|^{\frac{m}{n}} \right \rceil
\\&\leq
\sum_{
\left(c_1,c_3\right)\in \left\{-H,\ldots,H\right\}^2
\atop
r_{m,n}\left|c_1\right|^{\frac{n-m}{n}}\left|c_3\right|^{\frac{m}{n}}\leq H}
2\left \lceil \frac{H}{\log \left(dH\right)} \right \rceil
\leq
8H^3\left(\frac{1}{\log \left(dH\right)}+\frac{1}{H}\right)
\end{align*}
where $\left \lceil x \right \rceil$ is the smallest integer no less than $x$. Therefore, the fraction of trinomials that are ill-conditioned is at most $\frac{1}{\log \left(dH\right)}+\frac{1}{H}$.

\section{Signs and Solving are Roughly Equivalent for Trinomials} 
\label{sec:equiv} 
In what follows, it clearly suffices to focus on sign evaluation on $\R_+$ and approximation of roots in $\R_+$, since we can simply work with $f(\pm x)$. 

\medskip 
\noindent 
{\bf Proof of Lemma \ref{lemma:equiv}:} ($\Longrightarrow$) If Koiran's Trinomial Sign Problem can be solved in polynomial-time then we simply apply bisection to solve for all the positive roots of any input trinomial $f\!\in\!\Z[x]$ with degree $d$ and all coefficients having logarithmic height $\log H$. 

In particular, Lemma \ref{lemma:cauchy} tells us that the positive roots of $f$ lie in the interval $(1/(2H),2H)$, and Lemma \ref{lemma:count} (combined with Baker's Theorem and Theorem \ref{thm:log}) tells us that 
we can count exactly how many positive roots there 
are in time $\log^{2+o(1)}(dH)$. 

If there are no positive roots then we are done. 

If there is only one positive root then we start with the interval $[0,2H]$ and then apply bisection (employing the assumed solution to Koiran's Trinomial Sign Problem) until we reach an interval of width 
$\frac{1}{2H\cdot 4(d-1)^2}$. Theorem \ref{thm:alpha} (combined with a simpler argument for the case $d\!=\!2$ involving completing the square and Theorem \ref{thm:binor}) then tells us that this is sufficient accuracy to obtain an approximate root in the sense of Smale. The number of bisection steps is clearly $O(\log(dH))$, so the overall final complexity is $\log^{O(1)}(dH)$ since we've assumed 
trinomial sign evaluation takes time $\log^{O(1)}(dH)$. 

If there are two positive roots then we first approximate the unique positive critical point $w$ of $f$ via Theorem \ref{thm:binor} (since it is the unique root of a binomial with coefficients of logarithmic height $O(\log(dH))$), to accuracy $\frac{1}{3\cdot 2H\cdot 
4(d-1)^2}$. This ensures that the intervals 
$(0,w)$ and $(w,2H)$ each contain a positive root of 
$f$, thanks to Theorem \ref{thm:alpha}. We then apply bisection (in each interval) as in the case of just one positive root, clearly ending in time $\log^{O(1)}(dH)$. \qed 

\medskip 
\noindent 
{\bf ($\Longleftarrow$):} Our argument is almost identical to the converse case, except that we 
re-organize our work slightly differently. First, 
we count the number of positive roots of $f$ in 
time $\log^{2+o(1)}(dH)$, as outlined in the 
converse case. Then, we additionally compute the 
signs of $f(0)$ and $f(+\infty)$ essentially for free by simply evaluation $\sign(c_1)$ and $\sign(c_3)$, 

This data will partition $\R_+$ into at most 
$3$ open intervals upon which $f$ has constant (nonzero) sign. 
So to evaluate $\sign(f(r))$ at an $r\!\in\!\Q$ with 
logarithmic height $\log H'$, we simply need to check which interval contains $r$ {\em or} if $r$ is itself a rational root of $f$. 

Doing the latter is already known algorithmically, thanks to earlier work of Lenstra: \cite{lenstra} in fact details a polynomial time algorithm for finding all rational roots of any sparse polynomial (and even extends to finding all bounded degree factors, over number fields of bounded degree). In particular, \cite{lenstra} also proves that the logarithmic heights of the rational roots of $f$ are $\log^{O(1)}(dH)$. So we can simply compare $r$ to the rational roots of $f$ (in time $\log^{O(1)}(dH)$) and output $\sign(f(r))\!=\!0$ if $r$ matches any such root. 

So let us now assume that $r$ is not a root of $f$. If we compute the positive roots of $f$ to accuracy $\frac{1}{3H'}$, then we can easily decide which interval contains $r$ and immediately compute $\sign(f(r))$. By assumption, finding a set of positive approximate roots respectively converging to each true positive root of $f$ is doable in time $\log^{O(1)}(dH)$. So, to potentially upgrade our approximate roots to accuracy $\frac{1}{3H'}$, we merely apply Newton Iteration: This involves $O(\log\log(dHH'))$ further iterations, and each such iteration involves $\log(dH)$ arithmetic operations. Since we only need accuracy $\frac{1}{3H'}$, we can in fact work with just the $O(\log(dHH'))$ most significant digits of our approximate roots. So we are done. \qed

\medskip 
\noindent 
{\bf Proof of Corollary \ref{cor:koiran}:} Our corollary follows easily from the proof of Lemma 
\ref{lemma:equiv}: Just as in the proof of 
the ($\Longleftarrow$) direction of our last proof, 
we first apply \cite{lenstra} to check whether  
$u/v$ is a root of $f$ in time $\log^{O(1)}(dH)$. 
If so, then we are done (with $\sign(f(u/v))\!=\!0$), so let us assume $u/v$ is {\em not} a root of $f$. 

From Theorem \ref{thm:big}, if $f$ is {\em not} ill-conditioned, then we can find a set of positive 
approximate roots respectively converging to each 
true positive root of $f$. Furthermore, by construction (by Algorithm \ref{algor:trinosolr} in particular), each approximate root is within distance $\frac{\zeta}{4(d-1)^2}$ of a (unique) true positive root $\zeta$. So, as before, we merely need to refine slightly via Newton iteration until we attain accuracy $\frac{1}{3H}$, to find which interval (on which $f$ has constant sign) $r$ lies in. This additional work takes time $\log^{4+o(1)}(dH)$, so 
we are done. \qed 

\bibliographystyle{plain}
\bibliography{20250501arxiv}

\end{document}